\documentclass[12pt]{amsart} 
\usepackage{amssymb}
\usepackage{latexsym}
\usepackage{amsmath}
\usepackage{amsthm}
\usepackage{amsfonts}
\usepackage{mathrsfs} 
\usepackage{soul}
\usepackage[usenames,dvipsnames]{xcolor}
\usepackage{doi}
\usepackage{tikz-cd} 
\usepackage{lineno} 
\usepackage{datetime} 
\usepackage{ulem}  
\usepackage{cancel} 

\newtheorem{theorem}{Theorem}[section]

\newtheorem{lemma}[theorem]{Lemma}
\newtheorem{corollary}[theorem]{Corollary}
\newtheorem{proposition}[theorem]{Proposition}

\theoremstyle{definition}

\theoremstyle{remark}
\newtheorem{remark}[theorem]{Remark}

\numberwithin{equation}{section}
\numberwithin{theorem}{section}
\allowdisplaybreaks
\usepackage{color}

\newcommand{\R}{\mathbb R}
\newcommand{\N}{\mathbb N}
\newcommand{\C}{\mathbb C}

\newcommand{\logl}{L\textnormal{log} L}

\newcommand{\lexp}{L_{\textnormal{exp}}}

\newcommand{\tlog}{T_{\textnormal{log}}}

\newcommand{\inter}{(-1,1)}
\newcommand{\pint}{(-1,1)}
\newcommand{\simb}{\text{sim }\mathcal{B}}

\newcommand{\mcB}{\mathcal B}
\newcommand{\intt}{ \int_{-1}^1}
\newcommand{\mlog}{m_{L^1}}
\newcommand{\lmlog}{L^1(\mlog)}
\newcommand{\nin}{{n=1}^\infty}
\newcommand{\kin}{{k=1}^\infty}
\newcommand{\hT}{\widehat T}
\newcommand{\mexp}{m_{\textnormal{exp}}}

\title[Measure theoretic aspects of the finite   Hilbert transform]
{Measure theoretic aspects of the \\ finite   Hilbert transform}

\author[G. P.  Curbera]{Guillermo P. Curbera}
\address{Facultad de Matem\'aticas \& IMUS,
Universidad de Sevilla, 
Calle Tarfia s/n,  Sevilla 41012, Spain}
\email{curbera@us.es}

\author[S. Okada]{Susumu Okada}
\address{112 Marcorni Crescent, Kambah, ACT 2902, Australia}
\email{sus.okada@outlook.com}

\author[W. J. Ricker]{Werner J. Ricker}
\address{Math.--Geogr.\  Fakult\"at, Katholische Universit\"at
Eichst\"att--Ingolstadt, D--85072 Eichst\"att, Germany}
\email{werner.ricker@ku.de}

\thanks{The first author acknowledges the support  of 
of PID2021-124332NB-C21 FEDER/Ministerio de Ciencia e Innovaci\'on and 
FQM-262 (Spain).}

\date{\today}


\begin{document}


\begin{abstract} 
The finite Hilbert transform $T$, when acting in the classical Zygmund space $\logl$
(over $(-1,1)$), was intensively studied in \cite{curbera-okada-ricker-log}. In this note
an integral representation of $T$ is established via the $L^1(-1,1)$-valued measure
$\mlog\colon A\mapsto T(\chi_A)$  for each Borel set $A\subseteq(-1,1)$. 
This integral representation,
together with various non-trivial properties of $\mlog$, allow the use of
measure theoretic methods (not available in \cite{curbera-okada-ricker-log}) to establish 
new properties of $T$. For instance, as an operator between Banach function spaces
$T$ is not order bounded, it is not completely continuous and neither is it  weakly compact.
An appropriate Parseval formula for $T$ plays a crucial role.
\end{abstract}

\maketitle


\section{Introduction}


The finite Hilbert transform $T(f)$ of $f\in L^1(-1,1)$ is the  principal value integral
\begin{equation*}
T(f)(t):=\lim_{\varepsilon\to0^+} \frac{1}{\pi}
\left(\int_{-1}^{t-\varepsilon}+\int_{t+\varepsilon}^1\right) \frac{f(x)}{x-t}\,dx ,
\end{equation*}
which exists for a.e.\ $t\in(-1,1)$ and is a measurable function.
The resulting linear operator $T\colon L^1\pint\to L^0\pint$ is 
continuous, as a consequence of   Kolmogorov's theorem,
\cite[Theorem III.4.9(b)]{bennett-sharpley}.  It follows from M. Riesz's 
theorem, \cite[Theorem III.4.9(a)]{bennett-sharpley}, that the operator $T$ 
maps $L^p\pint$ into itself boundedly whenever $1 < p < \infty$.  This 
result remains valid for  a certain class of 
rearrangement invariant (briefly, r.i.)  spaces. Indeed,  given a r.i.\  
space $X$, let $\underline{\alpha}_X $ and $\overline{\alpha}_X$ denote 
its  lower and upper  Boyd indices, respectively 
(see \cite[Definition III.5.12]{bennett-sharpley}).  
From  \cite[pp.170--171]{krein-petunin-semenov} it follows that  
$T$ maps $X$ boundedly into itself  if and only if $X$ admits 
non-trivial Boyd indices, that is, if $0<\underline{\alpha}_X\le \overline{\alpha}_X<1$.  
For example, $\underline{\alpha}_{L^p\pint}= \overline{\alpha}_{L^p\pint}=1/p  $ 
for $1\le p \le \infty$,  \cite[Theorem IV.4.6]{bennett-sharpley}.  
For such r.i.\  spaces $X$, the resulting linear operator $T_X\colon X\to X$ 
has been studied in \cite{curbera-okada-ricker-ampa},
\cite{curbera-okada-ricker-mh} and \cite{curbera-okada-ricker-am}.

In \cite{curbera-okada-ricker-log} the operator $T$ is investigated when acting
in the classical Zygmund space $\logl:=\logl\pint$. This is a r.i.\ space on $\pint$
which is close to $L^1\pint$ in the sense that it contains all r.i.\ spaces having
non-trivial Boyd indices: the Boyd indices of $\logl$ itself are trivial
(both equal 1), which implies that $T$   \textit{does not} map $\logl$
into itself. However, $T$ \textit{does} map $\logl$ continuously into the
strictly larger space $L^1\pint$, \cite[Theorem 2.1]{curbera-okada-ricker-log}.
It is also established in Theorem 5.6 of \cite{curbera-okada-ricker-log}
that $T$ cannot be extended to any r.i.\ domain space beyond $\logl$ whilst
still taking its values in $L^1\pint$, that is, the finite Hilbert transform
$T\colon\logl\to L^1\pint$, denoted by $\tlog$, is optimally defined.

The point of departure of this paper is to continue the study of $\tlog$ but, from a different
perspective, which then allows the use of other techniques. A crucial fact is that 
the finitely additive 
set function $\mlog\colon\mcB\to L^1\pint$ defined on the Borel $\sigma$-algebra
$\mcB$ of $\pint$ by
$\mlog(A):=T(\chi_A)$, for $A\in\mcB$,
actually turns out to be $\sigma$-additive, that is, $\mlog$ is a \textit{vector measure}.
This allows us to invoke relevant results from the well developed theory of
integration of scalar functions with respect to Banach-space-valued measures;
see, for example, \cite{diestel}, \cite{lewis}, \cite{ORR},  \cite{okada-ricker-sanchez} and the references 
therein.
In particular, there is available the associated Banach function space
$L^1(\mlog)$ of all ($\C$-valued) $\mlog$-integrable functions together
with the continuous linear \textit{integration operator} 
$I_{\mlog}\colon L^1(\mlog)\to L^1\pint$ given by 
\begin{equation}\label{1.1}
I_{\mlog}(f):=\int_{\pint}f\,d\mlog,\quad f\in L^1(\mlog).
\end{equation}
One of the main results is that $L^1(\mlog)$ turns out to coincide with $\logl$ and that
$I_{\mlog}(f)=\tlog(f)$ for every $f\in L^1(\mlog)$, that is, \eqref{1.1}
provides an \textit{integral representation} of the singular integral 
operator $\tlog$ (cf.\ Theorem \ref{t-4.1}).
This representation is used to establish various properties of $\tlog$.
The fact that $\mlog$ is a rather complicated vector measure plays an important role.
For instance, $\mlog$ has infinite variation over every set $A\in\mcB$ 
with positive Lebesgue
measure and the subset $\mlog(\{A\cap B:B\in\mcB\})$ is not 
order bounded in $L^1\pint$; see Theorem \ref{t-3.4} and 
Corollary \ref{c-3.6}, respectively.
Moreover, the range $\mlog(\mcB)$ of $\mlog$ is not a relatively 
compact subset of $L^1\pint$; see Theorem \ref{t-3.10}. Some 
consequences of such properties of $\mlog$,
together with \eqref{1.1} and the identity $I_{\mlog}=\tlog$, are 
as follows.
The finite Hilbert transform $\tlog\colon\logl\to L^1\pint$ fails 
to be an order bounded operator
(cf.\ Theorem \ref{t-4.4}), it is not a completely continuous 
operator (cf.\ Theorem \ref{t-4.6}) and
neither is it a weakly compact operator (cf.\ Theorem \ref{t-4.7}).


\section{Preliminaries}


\subsection{Rearrangement invariant spaces}


We  denote by $\mu$  the Lebesgue measure on $\R$. 
The relevant measure space considered in this paper is 
$\mu$ restricted to the Borel $\sigma$-algebra $\mcB$ 
of the open interval $(-1,1)$,
which we again denote by $\mu$. The vector space of all $\C$-valued, $\mcB$-simple functions 
is denoted by $\simb$.   We denote by $L^0:=L^0(\mu)$ the vector space of all 
$\C$-valued measurable functions on $(-1,1)$.  Measurable functions 
on $(-1,1)$ which  
coincide a.e.\  are identified.  With respect to the a.e.\  pointwise order
for its positive cone, $L^0$ 
is a complex vector lattice.  The space $L^0$ is also a metrizable topological
vector space for the topology of convergence in measure.  
An order ideal  $X$ of $L^0$ is called a \textit{Banach 
function space} (briefly, B.f.s.) based on the measure space $((-1,1), \mcB, \mu)$ 
if $\simb\subseteq X$ and if $X$ is equipped with a lattice norm $\|\cdot\|_X$ for which 
it is complete.  For brevity, we speak of  $X$ as a B.f.s.\ over $(-1,1)$.  A typical 
example of a B.f.s.\ over $ (-1,1)$ is the Lebesgue space $L^p:= L^p(\mu)$ for 
each $ 1 \le p \le \infty$.

A B.f.s.\ $X$ has \textit{$\sigma$-order continuous} norm if $\|f_n\|_X \downarrow 0$ 
for every sequence $(f_n)_{n=1}^\infty\subseteq X$ satisfying 
$f_n\downarrow 0$.  Equivalently, every function $f\in X$ has 
absolutely continuous norm, meaning that  $\|f\chi_{A_n}\| \downarrow 0$ 
whenever $A_n \downarrow\emptyset$ in $\mcB$ 
(see \cite[Propositions I.3.2 and I.3.5]{bennett-sharpley}).
We say that $X$ has the \textit{Fatou property} if, whenever $(f_n)_{n=1}^\infty$ 
is a norm-bounded, increasing sequence of non-negative functions in $X$ for which
$f=\sup_{n\in\N} f_n $ exists in $L^0$, then necessarily  $f \in X$ 
and $\|f_n\|_X\uparrow \|f\|_X$.
The \textit{associate space} $X'$ of $X$  is  the vector sublattice of $L^0$ 
defined by
$X':= \{g\in L^0: fg \in L^1 \text{  for every  } f\in X\}$
equipped with the lattice norm $\|g\|_{X'} : = \sup\{|\intt fg\,d\mu|: 
f \in X \text{ with   } \|f\|_X\le 1\}$.  It turns out that $X'$ is a closed linear subspace of the dual 
Banach space $X^*$ of $X$.  Moreover, $X'$ is also a B.f.s.\ over $(-1,1)$.  
If $X$ has  $\sigma$-order  continuous norm, then $X'= X^*$; 
see, for example,  \cite[p.29]{lindenstrauss-tzafriri}.

The decreasing rearrangement $f^*\colon [0,2]\to [0,\infty]$  
 of a function $f \in L^0\pint$  is the
right continuous inverse of its distribution function:
$ \lambda\mapsto \mu\big(\{t\in (-1,1):\,|f(t)|>\lambda\}\big)$ for $\lambda\ge 0$.
A \textit{rearrangement invariant}  space $X$ over $(-1,1)$ is a
B.f.s.\ with the Fatou property
such that if $g^*\le f^*$ with $f\in X$,
then $g\in X$ and $\|g\|_X\le\|f\|_X$. In this case, its associate space $X'$ 
is also a r.i.\  space, \cite[Proposition II.4.2]{bennett-sharpley}.

The Zygmund space $\logl$ is the vector sublattice of $L^0$ 
consisting of all measurable functions $f$
on $(-1,1)$ for which either one of the following two equivalent conditions holds:
$$
\int_{-1}^1|f(x)|\log^+|f(x)|\,dx<\infty,\quad \int_0^2 f^*(t)\log\Big(\frac{2e}{t}\Big)\,dt<\infty,
$$
see \cite[Definition IV.6.1 and Lemma IV.6.2]{bennett-sharpley}.
Then  $\logl$ is a  r.i.\  space   with $\sigma$-order  continuous norm given by
$$
\|f\|_{\logl}:=\int_0^2 f^*(t)\log\Big(\frac{2e}{t}\Big)\,dt,\quad f\in\loglº.
$$
(cf. \cite[Definition IV.6.3]{bennett-sharpley}).
The Zygmund space $\lexp\subseteq L^0$  
is defined as the  vector sublattice of all measurable functions $f$ 
on $(-1,1)$ for which there is a constant $\lambda> 0 $ (depending on $f$) such that
$$
\intt \exp\big(\lambda|f(x)|\big)\,dx< \infty;
$$
see \cite[Definition IV.6.1]{bennett-sharpley}.
It is a r.i.\  space when equipped with the lattice norm 
\begin{equation}\label{2.1}
\|f\|_{\lexp}: = \sup_{0<t<2}\frac{f^{**}(t)}{\log(2e/t)},\quad f\in \lexp.
\end{equation}
Here, $f^{**}(t):=(1/t)\int_0^tf^*(s)\,ds$, for $t>0$, is the maximal function of $f^*$; 
see \cite[Definition II.3.1]{bennett-sharpley}.
The space $ \lexp$ is the associate space of $\logl$ 
(up to equivalence of norms),   
\cite[Theorem  IV.6.5]{bennett-sharpley}.  Since $\logl$ has $\sigma$-order  
continuous norm, 
it follows that
\begin{equation}\label{2.2}
\lexp = (\logl)' = (\logl)^*.
\end{equation}
Moreover, the  inclusions 
$L^\infty \subseteq \lexp \subseteq L^p \subseteq 
\logl\subseteq L^1$ are valid,
for every $1< p < \infty$, and the corresponding natural inclusion maps are bounded; 
see again  \cite[Theorem  IV.6.5]{bennett-sharpley}.


\subsection{The finite Hilbert transform}


Since $0<\underline{\alpha}_{\logl}=  \overline{\alpha}_{\logl}= 1$ 
(see \cite[Theorem IV.6.5]{bennett-sharpley}), it is the case that
$T(\logl) \not\subseteq\logl$. Otherwise, the inclusion 
$T(\logl) \subseteq\logl$ would imply that $T$ maps 
$\logl$ boundedly into itself, which is not the case; see Section 1.  
However, the finite Hilbert transform $T$ does map $\logl$ boundedly into 
$L^1$  as ascertained in  \cite[Theorem 2.1]{curbera-okada-ricker-log}, 
which we record formally.

\begin{lemma}\label{l-2.1}
The finite Hilbert transform $T$ maps $\logl$ into $L^1$ 
and the resulting linear operator
$\tlog\colon  \logl \to L^1 $ is bounded.
\end{lemma}

The following Parseval formula for a pair of functions $f \in \logl$ and $g \in L^\infty$, 
established in \cite[Corollary  3.3]{curbera-okada-ricker-log}, will 
play an important role in the sequel. 
An analogous Parseval formula involving a r.i.\  
space $X$ with non-trivial Boyd indices and its associate space $X'$ 
occurs in \cite[Proposition 3.1]{curbera-okada-ricker-ampa}.


\begin{lemma}\label{l-2.2}  
Let $f\in \logl$ and $ g\in L^\infty$. Both of the functions
$fT(g)$ and $gT(f)$ belong to $L^1$ and satisfy
\begin{equation}\label{2.3}
\int_{-1}^1fT(g)\,d\mu=-\int_{-1}^1gT(f)\,d\mu.
\end{equation}
\end{lemma}

From Lemma \ref{l-2.2}, given $g \in L^\infty$, the function 
$fT(g) \in L^1$ for every $f \in \logl$.  
That is, $T(g) \in (\logl)' = \lexp$.  Consequently, $T$ maps $L^\infty$ into 
$\lexp$.  Let 
\begin{equation}\label{2.4}
T_\infty\colon  L^\infty \to \lexp
\end{equation}
denote the resulting linear operator, which is 
necessarily bounded via the closed graph theorem.  It follows from  \eqref {2.2}
and \eqref{2.3}  that the adjoint  operator 
$(\tlog)^*\colon  (L^1)^* \to (\logl)^* = \lexp $ of $\tlog$ 
equals $-T_{\infty}$ on $L^\infty=(L^1)^*$, that is,
\begin{equation}\label{2.5}
(\tlog)^*(f) = -T_{\infty}(f) \in \lexp, \quad f \in L^\infty.
\end{equation}


Recall that given $\lambda\in(0,1]$ a function $\phi\colon(-1,1)\to\C$ is called
\textit{$\lambda$-H\"older continuous} if there exists a constant $K_\phi>0$ such that
$$
|\phi(x)-\phi(t)|\le K_\phi |x-t|^\lambda,\quad x,t\in(-1,1).
$$
Let $B(\cdot,\cdot)$ denote the Beta function and define 
the bounded, continuous function $w(x):=\sqrt{1-x^2}$ for 
$x\in(-1,1)$. The arcsine distribution $1/w$ belongs to $L^p$
for all $1\le p<2$ and satisfies
\begin{equation}\label{2.6}
\Big(T\Big(\frac{1}{w}\Big)\Big)(t)=\text{p.v. } 
\frac1\pi \int_{-1}^1 \frac{dx}{\sqrt{1-x^2}\, (x-t)}=0,
\quad t\in(-1,1).
\end{equation}
\cite[\S4.3(14)]{tricomi}. The following result, together with its proof, 
was kindly communicated to us (privately) by Prof. D. Elliott.

\begin{theorem}\label{t-2.3}
Let $\lambda\in(0,1]$ and $\phi\colon(-1,1)\to\C$ be a 
$\lambda$-H\"older continuous function. Then
\begin{equation}\label{2.7}
w(t)\Big|\Big(T\Big(\frac{\phi}{w}\Big)\Big)(t)\Big|\le \frac{2}{\pi} K_\phi\, B(1/2,\lambda),
\quad t\in(-1,1).
\end{equation}
In particular, $wT(\phi/w)\in L^\infty$.
\end{theorem}

\begin{proof}
Fix $t\in(-1,1)$. The proof is via a series of steps.

\textit{Step 1.} Define $I_1(t):=  \int_{-1}^t (t-x)^{\lambda-1}/\sqrt{1-x^2}\,dx$ and
note that the integrand is non-negative on $(-1,t)$. For the change of variables
$x=-1+(1+t)\xi$, that is, $\xi=(1+x)/(1+t)$, we have that 
 $(t-x)=(1+t)(1-\xi)$ and $(1-x^2)=2(1+t)\xi(1-(1+t)(\xi/2))$. It follows that
\begin{equation}\label{2.8}
I_1(t)=\frac{(1+t)^{\lambda-\frac12}}{\sqrt2}
\int_0^1\frac{(1-\xi)^{\lambda-1}}{\sqrt\xi\,\sqrt{1-(1+t)(\xi/2)}}d\xi.
\end{equation}
For each $\xi\in[0,1]$ the inequalities $0\le (1+t)(\xi/2)\le (1+t)/2$ imply that
$$
1- (1+t)(\xi/2)-(1-t)/2\ge 1-(1+t)/2-(1-t)/2=0.
$$
Accordingly,
\begin{equation}\label{2.9}
0<\Big(1-(1+t)(\xi/2)\Big)^{-1/2}\le \frac{\sqrt2}{\sqrt{1-t}},\quad \xi\in[0,1].
\end{equation}
It follows from \eqref{2.8} and \eqref{2.9} that
$$
I_1(t)\le \frac{(1+t)^{\lambda-\frac12}}{\sqrt2}\int_0^1\frac{\sqrt2\,(1-\xi)^{\lambda-1}}
{\sqrt\xi\,\sqrt{1-t}}\,d\xi=
\frac{(1+t)^{\lambda}}{\sqrt{1-t^2}}
\int_0^1\xi^{\frac12 -1}(1-\xi)^{\lambda-1}\,d\xi
$$
and hence, via the definition of the Beta function, that
\begin{equation}\label{2.10}
I_1(t)\le B(1/2,\lambda)\frac{(1+t)^\lambda}{\sqrt{1-t^2}}
=B(1/2,\lambda)\frac{(1+t)^\lambda}{w(t)}.
\end{equation}

\textit{Step 2.}  Define $I_2(t):=  \int_t^1 (x-t)^{\lambda-1}/\sqrt{1-x^2}\,dx$
and note  that the integrand is non-negative on $(t,1)$. Making the substitution
$s=-x$ in the integral 
$$
I_1(-t)=\int_{-1}^{-t}\frac{(-t-x)^{\lambda-1}}{\sqrt{1-x^2}}\,dx
=\int_{t}^{1}\frac{(-t+s)^{\lambda-1}}{\sqrt{1-s^2}}\,ds
$$
shows that $I_1(-t)=I_2(t)$. Since $w(t)=w(-t)$, it follows from the previous identity 
and \eqref{2.10} that
\begin{equation}\label{2.11}
I_2(t)=I_1(-t)\le B(1/2,\lambda)\frac{(1-t)^\lambda}{w(t)}.
\end{equation}

\textit{Step 3.}  In view of \eqref{2.6} we can write
\begin{equation}\label{2.12}
\Big(T\Big(\frac{\phi}{w}\Big)\Big)(t)=\text{p.v. } \frac1\pi \int_{-1}^1 \frac{\phi(x)-\phi(t)}{\sqrt{1-x^2}\,(x-t)}\,dx.
\end{equation}
Moreover, since $|x-t|$ equals $(t-x)$ for $-1<x\le t$ and equals
$(x-t)$ for $t\le x<1$, the $\lambda$-H\"older continuity of $\phi$ together with
\eqref{2.10} and \eqref{2.11} yields
\begin{align*}
\frac1\pi \int_{-1}^1 \Big| \frac{\phi(x)-\phi(t)}{\sqrt{1-x^2}\,(x-t)}\Big| \,dx
& \le 
\frac{K_\phi}{\pi} \int_{-1}^1 \frac{|x-t|^\lambda}{\sqrt{1-x^2}\,|x-t|}\,dx
\\ & =
\frac{K_\phi}{\pi} (I_1(t)+I_2(t))
\\ & \le
\frac{K_\phi}{\pi} B(1/2,\lambda)\frac{h(t)}{w(t)},
\end{align*}
where $h(s):=(1+s)^\lambda +(1-s)^\lambda$ for $s\in(-1,1)$.
Defining $h(-1)=h(1)=2^\lambda$ it is clear that $h\ge0$ is continuous
on $[-1,1]$, differentiable on $(-1,1)$ and satisfies $h(s)=h(-s)$. 
Since $h'(s)<0$ for $s\in(0,1)$, we can conclude that $h$ is decreasing 
on $[0,1]$ and so 
$
\sup_{|s|<1} h(s)=\max_{x\in[0,1]} h(x)=h(0)=2.
$
It follows that
$$
\frac1\pi \int_{-1}^1 \Big| \frac{\phi(x)-\phi(t)}{\sqrt{1-x^2}\,(x-t)}\Big| \,dx
\le  \frac{2K_\phi B(1/2,\lambda) }{\pi w(t)}.
$$
Accordingly, the principal value integral in \eqref{2.12} is a Lebesgue 
integral and so the previous inequality yields the desired inequality \eqref{2.7}.
\end{proof}


\section{The vector measure associated to $T$}


By a \textit{vector measure} we mean  a 
Banach-space-valued,  $\sigma$-additive set function defined on a $\sigma$-algebra 
of  subsets of some non-empty set.   
Given any  r.i.\ space $X$ over $\pint$ with non-trivial Boyd indices, 
the finite Hilbert transform $T_X\colon X\to X$ is bounded. 
The associated set function $A\mapsto m_X(A):=T_X(\chi_A)\in X$, for 
$A\in\mcB$, is clearly finitely additive. 
It turns out that $m_X$ is actually $\sigma$-additive, that is,
$m_X$ is a vector measure; see \cite[Section 3]{curbera-okada-ricker-qm}.


Since $T(L^\infty)\subseteq T(L^2)\subseteq L^1$ we can consider  the $L^1$-valued 
set function defined by
\begin{equation}\label{3.1}
\mlog(A): = T(\chi_A) \in L^1, \quad A \in \mcB.
\end{equation}
Clearly, $\mlog$ is finitely additive. Actually $\mlog$ is $\sigma$-additive.
Indeed, if $A_n\downarrow\emptyset$ in $\mcB$, then $\chi_{A_n}\to0$  
pointwise in $\pint$. The $\sigma$-order continuity of the  norm in $L^2$ implies that
$\chi_{A_n}\to0$ in $L^2$. Since $T_{L^2}$ is continuous, it follows 
that  $\mlog(A_n)=T(\chi_{A_n})\to0$ in $L^2$ and hence, also in $L^1$. 
So, $\mlog\colon  \mcB\to L^1$  is a vector measure.   
A subset $A\in\mcB$ is called $\mlog$-\textit{null } if $\mlog(B)=0$ 
for all sets $B\in \mcB$ with $B \subseteq A$.
Given $g \in L^\infty= (L^1)^*$, define a $\C$-valued  measure
$\langle \mlog,\,g\rangle$ on $\mcB$ by
\begin{equation*}
\langle \mlog,\,g\rangle(A): = \langle \mlog(A),\,g\rangle,\quad A \in  \mcB.
\end{equation*}


\begin{proposition}\label{p-3.1}
The following statements hold for the vector measure $\mlog\colon \mcB\to L^1$.
\begin{itemize}
\item[(i)]  For each $ g \in L^\infty$ it is the case that
\begin{equation}\label{3.2}
\langle \mlog, g\rangle(A) = -\int_AT(g)\,d\mu,\quad A \in {\mathcal B}.
\end{equation}
In particular, the variation measure $|\langle \mlog, g\rangle |$ of
$\langle \mlog, g\rangle$ is given by
\begin{equation}\label{3.3}
|\langle \mlog, g\rangle |(A) = \int_A|T(g)|\,d\mu,\quad A \in {\mathcal B}.
\end{equation}
\item[(ii)] The vector measure $\mlog$ and Lebesgue measure $\mu$ 
have the same null sets.
\end{itemize}
\end{proposition}

\begin{proof}
(i) Fix $g\in L^\infty$. Apply Lemma \ref{l-2.2} with $f:= \chi_A$  
to obtain \eqref{3.2} as
$$
\langle \mlog,g\rangle (A) = \langle \mlog(A), g\rangle=  
\int_{-1}^1 T(\chi_A)g\,d\mu = -  \int_{-1}^1 \chi_AT(g)\,d\mu = -\int_AT(g)\,d\mu,
$$
for every $A \in {\mathcal B}$.
It is well known that \eqref{3.2} implies \eqref{3.3}.

(ii) It is obvious from \eqref{3.1}  that every $\mu$-null set is $\mlog$-null. 
Conversely fix an $\mlog$-null set $A \in {\mathcal B}$. Then 
$|\langle \mlog,\, g\rangle|(A) = 0 $ for every $g \in L^\infty$. 
For $g: =\chi_{(-1,1)} \in L^\infty$, it follows from \eqref{3.3} that
$$
0 = |\langle \mlog,\chi_{(-1,1)} \rangle|(A) =\int_A |T(\chi_{(-1,1)})| d\mu 
= \frac{1}{\pi}\int_A\left|\,\log\left(\frac{1-t}{1+t}\right)\right|dt.
$$
This implies that $\mu(A) =0$ because $\log\left(\frac{1-t}{1+t}\right)=
0$ if and only if $t=0$.
\end{proof}


For the vector measure $m_X$, with $X$  a r.i.\ space 
over $\pint$ having non-trivial Boyd indices, 
Proposition 3.2(iii) in \cite{curbera-okada-ricker-qm} 
reveals that there exists $g_0\in X'\subseteq X^*$ such that
$\mu=|\langle m_X,g_0\rangle|$. Such a function $g_0$ is called a 
\textit{Rybakov functional} for $m_X$, \cite[Theorem IX.2.2]{diestel},
because $m_X$ and $|\langle m_X,g_0\rangle|$ have the same null sets.
The following result shows that an analogue is available for $X=L^1$.


\begin{proposition}\label{p-3.2}
 There exists a Rybakov functional $g_0\in L^\infty$ for
$\mlog$ satisfying $\mu=|\langle\mlog,g_0\rangle|$.
\end{proposition}

\begin{proof}
The proof is via a series of steps.

\textit{Step 1.} Define the function $\Psi\colon(0,1)\to\R$ by 
$$
\Psi(t):=\int_0^t\frac{1-(x/t)}{\sqrt{1-x^2}}\,dx+
\int_t^1\frac{(t-x)/(1-t)}{\sqrt{1-x^2}}\,dx,\quad t\in(0,1).
$$
The claim is that there exists a unique number $a\in(\frac12,1)$ such that $\Psi(a)=0$.
Since $x\mapsto 1/\sqrt{1-x^2}=1/w(x)$ belongs to $L^1(0,1)$ it follows
by  dominated convergence that $\Psi$ is continuous on $(0,1)$.

Given $t\in(0,1)$, the inequalities
$$
0\le \int_0^t\frac{1-(x/t)}{\sqrt{1-x^2}}\,dx
\le  \int_0^t\frac{1}{\sqrt{1-x^2}}\,dx
$$
imply that $\lim_{t\to0^+}\int_0^t\frac{1-(x/t)}{\sqrt{1-x^2}}\,dx=0$. Moreover,
for $t\to0^+$ we also have
$$
\int_t^1\frac{(t-x)/(1-t)}{\sqrt{1-x^2}}\,dx=
\frac{1}{1-t}\int_0^1\frac{(t-x)\chi_{[t,1]}(x)}{\sqrt{1-x^2}}\,dx
\to\int_0^1\frac{-x}{\sqrt{1-x^2}}\,dx<0,
$$
again by dominated convergence. Accordingly,
$$
\Psi(0^+)=\lim_{t\to0^+}\Psi(t)=-\int_0^1\frac{x}{\sqrt{1-x^2}}\,dx=-1<0.
$$

Concerning $\Psi(1^-)$ observe, again by dominated convergence, that
$$
\int_0^t\frac{1-(x/t)}{\sqrt{1-x^2}}\,dx=
\int_0^1\frac{1-(x/t)}{\sqrt{1-x^2}}\chi_{[0,t]}(x)\,dx
\to\int_0^1\frac{1-x}{\sqrt{1-x^2}}\,dx
$$
for $t\to1^-$. On the other hand, the inequality
$$
\Big|\frac{t-x}{1-t}\Big|\chi_{[t,1]}(x)=\frac{x-t}{1-t}\le \frac{1-t}{1-t}=1,\quad x\in[t,1],
$$
for each $t\in(0,1)$, implies that
$$
\Big|\frac{t-x}{(1-t)\sqrt{1-x^2}}\Big|\chi_{[t,1]}(x)
\le \frac{1}{\sqrt{1-x^2}},\quad x\in(0,1).
$$
Since $1/w\in L^1(0,1)$, we can again apply dominated convergence to conclude that
$\int_t^1\frac{t-x}{(1-t)\sqrt{1-x^2}}\,dx\to0$ for $t\to1^-$. Accordingly,
$$
\Psi(1^-)=\lim_{t\to1^-}\Psi(t)=\int_0^1\frac{1-x}{\sqrt{1-x^2}}\,dx>0.
$$
By the intermediate value theorem there exists $a\in(0,1)$ satisfying $\Psi(a)=0$.

Fix $0<t<s<1$. From the definition of $\Psi$ it follows that
$$
\Psi(s)-\Psi(t)=\int_0^t \frac{\frac{x}{t}-\frac{x}{s}}{w(x)}\,dx+
\int_t^s \frac{(1-\frac{x}{s})+\frac{x-t}{1-t}}{w(x)}\,dx
+\int_s^1\frac{(1-x)(s-t)}{(1-s)(1-t)w(x)}\,dx,
$$
from which it is clear that $\Psi(s)>\Psi(t)$. Hence, $\Psi$ is strictly increasing.
Since
$$
\Psi\Big(\frac12\Big)=\int_0^1\frac{1-2x}{\sqrt{1-x^2}}\,dx=
\int_0^1\frac{1}{\sqrt{1-x^2}}\,dx-\int_0^1\frac{2x}{\sqrt{1-x^2}}\,dx
=\frac{\pi}{2}-2<0,
$$
the unique number $a\in(0,1)$ such that $\Psi(a)=0$ satisfies $\frac12<a<1$.

\textit{Step 2.} For $a$ as in Step 1 define $c_1:=-a$, $c_2:=0$ and $c_3:=a$ and
then the disjoint intervals $I_1:=(-1,c_1]$, $I_2:=(c_1,c_2]$, $I_3:=(c_2,c_3]$ and
$I_4:=(c_3,1)$. Now define the piecewise affine function $u\colon(-1,1)\to\R$ by
$$
u(x):=\begin{cases} 
\frac{a+x}{1-a}, \text{ for $x\in I_1$,}\\
\frac{a+x}{a}, \text{ for $x\in I_2$,}\\
\frac{a-x}{a}, \text{ for $x\in I_3$,}\\
\frac{a-x}{1-a}, \text{ for $x\in I_4$.}\\
\end{cases}
$$
Setting $u(-1):=u(1):=-1$, it is clear that $u$ is an even function which
is continuous on $[-1,1]$ and satisfies $|u(x)|\le1$ for $x\in[-1,1]$. Via Step 1 we note that
\begin{align}\label{3.4}
\int_{-1}^1\frac{u(x)}{\sqrt{1-x^2}}\,dx&=2\int_{0}^1\frac{u(x)}{\sqrt{1-x^2}}\,dx
\\ & \nonumber =2\Big(\int_{0}^a\frac{1-(x/a)}{\sqrt{1-x^2}}\,dx+
\int_{a}^1\frac{(a-x)/(1-a)}{\sqrt{1-x^2}}\,dx\Big)
\\ & \nonumber  =2\Psi(a)=0.
\end{align}

The claim is that $u$ is 1-H\"older continuous on $(-1,1)$.
Since $u$ is
piecewise affine, it is differentiable in $\pint$ except at the points
$c_1, c_2, c_3$, with constant derivative  (at most $1/(1-a)$ in modulus) in
the interior of each interval $I_j$, for $j\in\{1,2,3,4\}$. Hence, 
given two points $x<t$ in $\pint$, an application of the 
mean value theorem yields
\begin{equation}\label{3.5}
|u(x)-u(t)|\le 
\frac{1}{1-a}|x-t|,\quad x,t\in(-1,1),
\end{equation}
by considering the possible cases that arise. Namely, both points $x,t\in I_j$ for some $1\le j\le 4$, or $x,t$ lie in adjacent intervals (that is, $x\in I_j$ and $t\in I_{j+1}$, for $j=1,2,3$),
or $x\in I_j$ and $t\in I_{j+2}$ for $j=1,2$, or $x\in I_1$ and $t\in I_4$. The 
$1$-H\"older continuity of $u$ is clear from \eqref{3.5}.

\textit{Step 3.} The claim is that $u$ satisfies the inequality
\begin{equation}\label{3.6}
\Big|\sqrt{1-u(x)^2}-\sqrt{1-u(t)^2}\Big|\le \Big(\frac{2}{1-a}\Big)^{1/2}|x-t|^{1/2},\quad
x,t\in (-1,1).
\end{equation}
To establish \eqref{3.6} will require the fact that
\begin{equation}\label{3.7}
|\sqrt y-\sqrt z|\le \sqrt{|y-z|},\quad y,z\in(0,\infty),
\end{equation}
which follows, for  $y>z$, by squaring both sides.
The inequality \eqref{3.7}
implies that the left-side of \eqref{3.6} is at most
\begin{align*}
\Big|\sqrt{1-u(x)^2}-\sqrt{1-u(t)^2}\Big|&\le
\sqrt{|u(x)^2-u(t)^2|}=\sqrt{|u(x)-u(t)|\cdot|u(x)+u(t)|}
\\ & \le \sqrt2 \sqrt{|u(x)-u(t)|},
\end{align*}
because $|u(s)|\le1$ for $s\in(-1,1)$; see Step 2. But,
$|u(x)-u(t)|\le \frac{1}{1-a}|x-t|$ (see \eqref{3.5}), which then implies the validity of \eqref{3.6}.

\textit{Step 4.} Define the function $F\colon(-1,1)\to\C$ by
\begin{equation}\label{3.8}
F(x):=\begin{cases} 
u(x)-i\sqrt{1-u(x)^2}, \text{ for $x\in (-1,0]$,}\\
u(x)+i\sqrt{1-u(x)^2}, \text{ for $x\in (0,1)$.}\\
\end{cases}
\end{equation}
Recall that $u\in C([-1,1])$. 
Moreover, $\lim_{x\to0^+}F(x)=\lim_{x\to0^-}F(x)=u(0)=1$ and so there exists
an extension $F\in C([-1,1])\subseteq L^\infty$. 
The claim is that $F$ is $\frac12$-H\"older continuous on $(-1,1)$ and satisfies
\begin{equation}\label{3.9}
\int_{-1}^1\frac{F(x)}{w(x)}\,dx=0.
\end{equation}

It follows from \eqref{3.5} that
$$
|u(x)-u(t)|\le \frac{|x-t|}{1-a}=\frac{|x-t|^{1/2}|x-t|^{1/2}}{1-a}
\le \frac{\sqrt2}{1-a}|x-t|^{1/2},\quad x,t\in(-1,1),
$$
that is, $\text{Re}(F)=u$ is  $\frac12$-H\"older continuous. So, it suffices to show
that $G:=\text{Im}(F)$ is also $\frac12$-H\"older continuous. To this end, consider 
the following cases.

If $t,x\in(-1,0]$ or $t,x\in[0,1)$, then  the inequality \eqref{3.6} shows that
$|G(x)-G(t)|\le (2/(1-a))^{1/2}|x-t|^{1/2}$.

Let $t\in(-1,0)$ and $x\in[0,1)$. From $u(0)=1$ and $|u|\le1$ on $(-1,1)$ it follows that
\begin{align*}
|G(x)-G(t)|&=\sqrt{1-u(x)^2}+\sqrt{1-u(t)^2}
\\ & =
\sqrt{u(0)-u(x)}\sqrt{u(0)+u(x)}+\sqrt{u(0)-u(t)}\sqrt{u(0)+u(t)}
\\ &\le 
\sqrt2 (\sqrt{u(0)-u(x)}+\sqrt{u(0)-u(t)}).
\end{align*}
Applying \eqref{3.5} and the inequalities $|x|\le |x-t|$ and $|t|\le |x-t|$ we have that 
\begin{equation}\label{3.10}
|G(x)-G(t)|\le \Big(\frac{2}{1-a}\Big)^{1/2}(|x|^{1/2}+|t|^{1/2})
\le 2 \Big(\frac{2}{1-a}\Big)^{1/2} |x-t|^{1/2}.
\end{equation}

In the case when  $t\in(0,1)$ and $x\in(-1,0)$ a similar 
argument as in the previous case again yields the inequality \eqref{3.10}.

We can conclude that \eqref{3.10} holds for all $x,t\in(-1,1)$. Accordingly, 
$G$ is $\frac12$-H\"older continuous on $(-1,1)$ and hence, so is $F=u+iG$.

To establish \eqref{3.9}, observe first that   $G\in C([-1,1])$ is an odd function
 and $1/w\in L^1$ is an even function,
so that $\int_{-1}^1 \frac{G(x)}{w(x)}\,dx=0$. This together with \eqref{3.4}
establish \eqref{3.9} because $F=u+iG$.

\textit{Step 5.} Setting $\phi:=F$ and $\lambda:=1/2$ in Theorem \ref{t-2.3} yields that
$g_0=-wT(F/w)\in L^\infty$. Fix $p\in(2,\infty)$. Since $F\in L^\infty$, we have
$F\in X:=L^p$ and, by \eqref{3.9}, that $\int_{-1}^1 \frac{F(x)}{w(x)}\,dx=0$. 
Also $g_0\in X$. Since $0<\underline{\alpha}_X=\overline{\alpha}_X<1/2$, it follows 
from \cite[Corollary 3.5(b)]{curbera-okada-ricker-ampa} with $f:=g_0$
and $g:=F$ there, that $g_0$ is the unique solution (in $X$) satisfying $T(g_0)=F$.

A direct calculation using \eqref{3.8} shows that $|F(x)|=1$ for all $x\in[-1,1]$. Accordingly,
$$
\mu(A)=\int_A|F|\,d\mu=\int_A|T(g_0)|\,d\mu,\quad A\in\mcB.
$$
On the other hand, \eqref{3.3}
implies that
$|\langle \mlog,g_0\rangle|(A)=\int_A|T(g_0)|\,d\mu$, for each $A\in\mcB$. 
So, $\mu=|\langle \mlog,g_0\rangle|$
with $g_0\in L^\infty=(L^1)^*$. Hence, $g_0$ is a Rybakov functional for $\mlog$
satisfying $\mu=|\langle \mlog,g_0\rangle|$.
\end{proof}


\begin{remark}\label{r-3.3}
In Proposition \ref{p-3.2} the proof of the existence of $g_0\in L^\infty$ 
which satisfies $\mu=|\langle\mlog,g_0\rangle|$  is significantly 
more involved than proving the existence of $g_0\in X'\subseteq X^*$,
which satisfies $\mu=|\langle m_X,g_0\rangle|$, whenever 
$X$ is a  r.i.\ space over $(-1,1)$ with non-trivial Boyd indices 
(cf. \cite[Proposition 3.2(iii)]{curbera-okada-ricker-qm}). On the other hand,
since $L^\infty\subseteq X'\subseteq X^*$ for all such r.i.\ spaces $X$ and
$m_X(A)=\mlog(A)$ for all $A\in\mcB$, Proposition \ref{p-3.2}
is a considerable strengthening
of Proposition 3.2(iii) in  \cite{curbera-okada-ricker-qm}. First,
$g_0$ exists in the \textit{proper subspace} $L^\infty$ of $X'$ 
and second, the \textit{same function} $g_0\in L^\infty$ 
can be chosen as a Rybakov functional satisfying $\mu=|\langle m_X,g_0\rangle|$
for \textit{every} $X$.
\end{remark}


The variation measure $|\mlog|\colon \mcB \to [0,\infty]$ of $\mlog$ is 
defined as for scalar measures by replacing the absolute value with the norm in $L^1$; 
see \cite[Definition I.1.4]{diestel}.


\begin{theorem}\label{t-3.4}
The vector measure  $\mlog \colon  \mcB \to L^1$
has infinite variation over every set
$A\in \mcB$ satisfying $\mu(A)>0$.
\end{theorem}

\begin{proof}  
Fix a set $A\in\mcB$ with  $\mu(A)>0$.  Then either $\mu(A\cap [0,1)) > 0$ 
or  $\mu(A\cap (-1,0)) > 0$ or both. Assume that $\mu(A\cap [0,1)) > 0$.  Fix $n \in \N$ with 
$n \ge 2$.  Given $k \in \N$ with $1 \le k < n$
we shall prove,
with $A_k(n):=A\cap[(k-1)/n,k/n)$, that 
\begin{equation}\label{3.11}
\|\mlog(A_k(n))\|_{L^1}  \ge \frac1\pi\int_{A_k(n)} \log (1-y)\,dy
+ \frac1\pi (\log n )\cdot\mu(A_k(n)).
\end{equation}
Since $n$ and $k$ are fixed,  let 
$a:= (k-1)/n$ and 
$b:= k/n$ for ease  of notation. Then 
\begin{align*}
\big\|T\big(\chi_{A_k(n)}\big)\big\|_{L^1} &=\big\|\chi_{(-1,b)}\cdot 
T\big(\chi_{A_k(n)}\big)\big\|_{L^1} + \big\|\chi_{[b,1)}\cdot 
T\big(\chi_{A_k(n)}\big)\big\|_{L^1}
\\ & \ge\big\|\chi_{[b,1)}\cdot T\big(\chi_{A_k(n)}\big)\big\|_{L^1}=
\int_b^1\left|T\big(\chi_{A_k(n)}\big)(x)\right| dx
\\ &=
 \frac{1}{\pi} \int_b^1 \left| \mathrm{p.v.}\!\intt \frac{\chi_{A_k(n)}(y)}{y-x}dy \right| dx
 = \frac{1}{\pi} \int_b^1 \left(\intt \frac{\chi_{A_k(n)}(y)}{x-y}dy \right) dx
 \\ & = \frac{1}{\pi} \intt \chi_{A_k(n)}(y)\left(\int_b^1\frac{1}{x-y}dx\right) dy
 =  \frac{1}{\pi}\int_{A_k(n)}\Big[\log(x-y)\Big]_{x=b}^{x=1}dy
\\ & =  \frac{1}{\pi} \int_{A_k(n)}\Big(\log(1-y)- \log(b-y)\Big)dy
\\ &\ge  
\frac{1}{\pi}\int_{A_k(n)}\Big(\log(1-y)- \log(b-a)\Big)dy
\\ & =   \frac{1}{\pi}\int_{A_k(n)} \log (1-y)\,dy -  
\frac{1}{\pi} \int_{A_k(n)}\log (1/n)\, dy
\\ &=  \frac{1}{\pi}\int_{A_k(n)} \log (1-y)\,dy 
\ + \  \frac{1}{\pi}(\log n)\cdot \mu(A_k(n)) ,
\end{align*}
where we  have used Fubini's theorem and the fact that $(x-y) \ge 0 $ 
whenever $b \le x < 1$ and $ y \in A_k(n)$. This establishes  
\eqref{3.11} because 
$\mlog(A_k(n)) =  T\big(\chi_{A_k(n)}\big)$ 
by definition.  Since, for every $n \in \N$ with $n\ge2$ 
the collection $\{A_k(n):1\le k<n\}$ is a partition 
of $A\cap[0,(n-1)/n)$ it follows that 
\begin{align*}
|\mlog|(A \cap [0,1))
& \ge 
|\mlog|\big(A \cap [0,(n-1)/n)\big)
\\ & \ge
\sum_{k=1}^{n-1}\left\|\mlog(A_k(n))\right\|_{L^1}
\\ & \ge 
\sum_{k=1}^{n-1} \left ( \frac{1}{\pi}\int_{A_k(n)} \log (1-y)\,dy \right.
+ \left.
\frac{1}{\pi}(\log n)\cdot \mu(A_k(n))\right) 
 \\ & = 
\frac{1}{\pi} \int_{A\cap ([0, (n-1)/n)} \log (1-y)\,dy 
+\frac{1}{\pi}(\log n)\cdot \mu\big(A\cap[0, (n-1)/n)\big).
\end{align*}
Since the function $y\mapsto \log(1-y)$ is integrable over $A\cap [0,1)$ 
and, for some  
$n_0\in\N$ with $n_0\ge2$ and $\alpha>0$, we have $\mu\big(A \cap [0,(n-1)/n)\big)>\alpha$ 
for all  $n\ge n_0$, it follows that $|\mlog|(A\cap[0,1)) = \infty$.

If $\mu(A\cap [0,1)) = 0$, then necessarily $\mu(A\cap (-1,0)) > 0$.   
By a similar argument, it follows that  $|\mlog|(A\cap(-1,0) )= \infty$.  
So, we can conclude that  $|\mlog|(A) = \infty$.
\end{proof}


\begin{remark}\label{r-3.5}
(i) Proposition \ref{p-3.1}(ii) shows that the vector measure
$\mlog$ is absolutely continuous with respect to $\mu$. Since $L^1$ does not have the 
Radon-Nikodym property \cite[p.219]{diestel}, the question arises 
of whether or not there exists a Bochner $\mu$-integrable (or, more general,
a Pettis $\mu$-integrable) function $F\colon\pint\to L^1$ such that 
$\mlog(A)=\int_AF\,d\mu$, for $A\in\mcB$. But,
the existence of such a function $F$ would imply that $\mlog$ has $\sigma$-finite
variation, \cite[Proposition 5.6(iv)]{vandulst}, which is not the case by Theorem
\ref{t-3.4}.

(ii) Given any r.i.\ space $X$ over $(-1,1)$ with non-trivial Boyd indices the 
operator $T_X$ is bounded on $X$. Since $X\subseteq L^1$ continuously,
\cite[Corollary II.6.7]{bennett-sharpley}, there exists $K>0$ such that
$\|g\|_{L^1}\le K\|g\|_X$ for all $g\in X$. It follows, for each $A\in\mcB$, that
$\|T(\chi_A)\|_{L^1}\le K \|T(\chi_A)\|_X$, that is, 
$\|\mlog(A)\|_{L^1}\le K \|m_X(A)\|_X$. So, Theorem \ref{t-3.4}
implies that $m_X$ has infinite variation over every set $A\in\mcB$ with $\mu(A)>0$.
For an alternative proof of this fact see \cite[Proposition 3.2(v)]{curbera-okada-ricker-qm}.
\end{remark}

Another  consequence of Theorem \ref{t-3.4} is as follows.  We will use the notation
$A\cap \mcB : = \{A\cap B: B \in \mcB\} =\{B: B \in \mcB,\;B \subseteq A\}$ for each 
$A \in \mcB$.
A subset $F\subseteq L^1$ is called \textit{order bounded} if there exists $0\le h\in L^1$
such that $|f|\le h$ for all $f\in F$.


\begin{corollary}\label{c-3.6} 
Given a set $A \in\mcB$ satisfying $\mu(A)>0$,  the subset $\mlog(A\cap \mcB )$ is not order bounded in $ L^1$. In particular, the range of $\mlog$ is not order bounded in $ L^1$.
\end{corollary}

\begin{proof}
Fix $A\in\mcB$ with $\mu(A)>0$.
Assume, by way of contradiction, that $\mlog(A\cap \mcB) $ is 
order bounded in $L^1$. Then the restriction of $\mlog$ to $A\cap\mcB$ 
is an $L^1$-valued order bounded measure defined on the measurable 
space $(A, A\cap\mcB)$.  Since $L^1$ is separable,  Dedekind 
complete and has $\sigma$-order continuous norm, it follows from  \cite[Proposition 1\,(1) and Theorem 5\,(2)]{SZ}
that there exists a smallest  vector measure $\nu\colon  A\cap\mcB \to L^1$,
taking values in the positive cone $(L^1)^+$ of $L^1$, which
dominates the restriction of $\mlog$ to   $A\cap\mcB$
in the sense that
$ |\mlog(B)| \le  \nu(B)$ in $(L^1)^+$ for each $ B \in A\cap  \mcB$. \ So, 
given $ B \in A\cap  \mcB$, we have $\|\mlog(B)\|_{L^1} \le \|\nu(B)\|_{L^1}$ 
because the $L^1$-norm is a lattice norm.  Hence, given finitely many 
pairwise disjoint sets $B_1,B_2,\dots, B_n \in A\cap  \mcB$ with $n \in \N$, 
it follows that
\begin{align*}
\sum_{k=1}^n \|\mlog(B_k)\|_{L^1} & 
\le
\sum_{k=1}^n \|\nu(B_k)\|_{L^1}
=\Big\|\sum_{k=1}^n \nu(B_k)\Big\|_{L^1}
\\ & =
\|\nu(\cup_{k=1}^n B_k)\|_{L^1}\le \|\nu(A)\|_{L^1} < \infty
\end{align*}
because $L^1$ is an abstract $L^1$-space, \cite[p.105]{okada-ricker-sanchez}, and because 
$\nu(B_k) \in(L^1)^+$  for $k=1,\dots, n$.  
Thus the total variation measure $|\mlog|\colon  \mcB \to [0,\infty]$ 
satisfies $|\mlog|(A) < \infty$.  This contradicts Theorem \ref{t-3.4}.  
The proof is thereby complete.

Applying the result just proven with $A=(-1,1)$, establishes the claim for the range of 
$\mlog$.
\end{proof}


\begin{remark}\label{r-3.7}
The results in the proof of Corollary \ref{c-3.6} that are quoted from
\cite{SZ} apply to the B.f.s.\ $L^1_\R$ of all $\mu$-integrable functions taking
their values in $\R$. Since $L^1=L^1_\R\oplus i L^1_\R$ is the complexification 
of $L^1_\R$, these results remain valid in $L^1$.
\end{remark}


To be able to establish the non-compactness of the
range $\mlog(\mcB)$ of $\mlog$  we need some preparation.

\begin{lemma}\label{l-3.8} 
The subset $U: = \{\chi_A: A \in \mcB\}$ of $L^0$ is not relatively compact.
\end{lemma}


\begin{proof}
Assume,  by way of contradiction, that $U$ is relatively compact in $ L^0$.  
This  will be shown to imply that $U$ is also relatively compact  in $L^1$.
Take any sequence $(A_n)_{n=1}^\infty$ in $\mcB$.  The relative compactness of 
$U$ in the metric space $L^0$ ensures that the corresponding sequence $(\chi_{A_n})_\nin$ 
admits a subsequence  $(\chi_{A_{n(k)}})_\kin$ converging in measure.   We may 
suppose that $(\chi_{A_{n(k)}})_\kin$ itself  converges pointwise a.e.\ on $(-1,1)$ 
by taking a subsequence, if necessary.  So, there exists $A\in \mcB$ such that 
$\chi_{A_{n(k)}}\to \chi_A$  pointwise a.e.\  as $k\to \infty$.   By dominated 
convergence, $\chi_{A_{n(k)}}\to \chi_A$  in $L^1$  as $k\to \infty$.   
In other words, an arbitrary sequence in the subset $U$ of $L^1$ admits a 
convergent subsequence in the norm of $L^1$.  This contradicts the fact 
that $U$ is not a relatively compact subset of $L^1$ 
(see, for example, \cite[Example III.1.2]{diestel}).  So, $U$ is not relatively compact in $L^0$.
\end{proof}


For the following fact we refer to \cite[Corollary 2.2.6]{bogachev}.

\begin{lemma}\label{l-3.9}  
For every function $h\in L^0$, the multiplication operator
$M_h\colon  f\mapsto fh$ is continuous from $L^0$ into itself.
\end{lemma}


Recall that
$w(x) :={\sqrt{1-x^2}}$ for $x \in (-1,1)$; see Section 2.
An important auxiliary operator $\widehat T$, given by
\begin{equation*}
\widehat{T}(g)(x):=-\frac{1}{w}\,
T(w\,g)(x),\quad \mathrm{a.e. }\; x\in (-1,1),
\end{equation*}
for each $g\in L^1$, will be required.  
It is a bounded linear operator from $L^p$ into itself whenever $1<  p < 2$, 
which is also  the case for every r.i.\  space $X$ over $\pint$ with Boyd indices 
satisfying  $1/2<\underline{\alpha}_X\le \overline{\alpha}_X<1$,
\cite[Theorem 3.2]{curbera-okada-ricker-ampa}.
Recall that  
$T\colon  L^1 \to L^0$ is continuous via   Kolmogorov's theorem.  Since $\hT$ 
is the composition of $T$ with the operator of multiplication by 
$w$ in $L^1$ and the operator of multiplication by $(-1/w)$ in $L^0$
(both continuous operators),  
it follows that $\hT$ maps 
$L^1$ continuously into $L^0$; we have used Lemma \ref{l-3.9}.


The Bartle-Dunford-Schwartz theorem implies that $\mlog(\mcB)$ is
necessarily a relatively weakly compact subset of $L^1$,
\cite[Corollary I.2.7]{diestel}. Corollary \ref{c-3.6} and the following 
result show that not more than this can be said about $\mlog(\mcB)$ as a subset of $L^1$.


\begin{theorem}\label{t-3.10}
The range $\mlog(\mcB)$ of $\mlog$ is not relatively compact in $L^1$.
\end{theorem}

\begin{proof}  Assume, on the contrary,  that  $\mlog(\mcB)$ is relatively 
compact in $L^1$.   Observe that $\mlog(\mcB)= T(U) \subseteq L^1$,
where $U$ is as in Lemma \ref{l-3.8}.  
As ascertained above,  $\hT$ maps $L^1$ continuously into $L^0$.  
So, $\hT (T(U))$  is relatively compact in $L^0$.   On the other hand, we 
know from \cite[Theorem 4.10 (iv)]{curbera-okada-ricker-log} that
\begin{equation}\nonumber
\chi_A - \frac{\mu(A)}{\pi w} = \hT (T(\chi_A)), \quad A \in \mcB,
\end{equation}
is valid as an equality in $\logl$ and hence, also in $L^0$.   It follows 
from the previous identity, for 
$V: = \{\chi_A- \frac{\mu(A)}{\pi w}: A \in \mcB\}$, 
that $V= \hT (T(U)) \subseteq L^0$. Accordingly, $V$ 
is relatively compact in $L^0$.  Moreover, the compact subset   
$W: =\{\frac{a}{\pi w}: 0 \le a \le 2\} \subseteq L^0$  satisfies
$U \subseteq V+W$.
This implies that $U$ is relatively  compact in $L^0$, which contradicts 
Lemma \ref{l-3.8}.   So, we can conclude that $\mlog(\mcB)$ is 
not relatively compact in $L^1$.
\end{proof}


Consider now the $\lexp$-valued 
set function associated with the operator $T_\infty\colon L^\infty\to \lexp$, defined by
\begin{equation*}
\mexp(A): = T_\infty(\chi_A) \in \lexp, \quad A \in \mcB.
\end{equation*}
Clearly, $\mexp$ is finitely additive and has bounded range. 
We will show that $\mexp$ is \textit{not} $\sigma$-additive.


\begin{theorem}\label{t-3.11} 
The bounded linear operator $T_\infty\colon L^\infty\to \lexp$  satisfies
\begin{equation*}
\|T_\infty(\chi_A)\|_{\lexp} > \frac{1}{e^2\pi}, \quad A \in \mcB,\,\, \mu(A)>0.
\end{equation*}
\end{theorem}

\begin{proof}
Let $f\in L^\infty$.  From  \eqref{2.1}  and
$f^{**}(t)=(1/t)\int_0^tf^*(s)\,ds\ge f^*(t)$, for $t>0$, it follows that 
\begin{equation}\label{3.12}
\|f\|_{\lexp}\ge\sup_{0<t<2}\frac{f^{*}(t)}{\log(2e/t)}.
\end{equation}
According to   \cite[Theorem 2.1]{edmunds-krbec} it is the case that
$$
\sup_{0<t<2}\frac{f^*(t)}{\log(2e/t)}<\infty \iff \sup_{1\le p<\infty} \frac{\|f\|_{L^p}}{p}<\infty
$$
and, from that proof, it follows for some constant $c>1/e$, that
\begin{equation}\label{3.13}
c \sup_{1\le p<\infty} \frac{\|f\|_{L^p}}{p}
\le \sup_{0<t<2}\frac{f^*(t)}{\log(2e/t)}\le e\, \sup_{1\le p<\infty} \frac{\|f\|_{L^p}}{p}.
\end{equation}
From \eqref{3.12} and \eqref{3.13} we can conclude that
\begin{equation}\label{3.14}
\|f\|_{\lexp} \ge c \, \sup_{1\le p<\infty} \frac{\|f\|_{L^p}}{p}.
\end{equation}

Let $E\subseteq\R$ be an arbitrary  measurable set with 
$\mu(E)<\infty$, $H$ be the Hilbert transform on $\R$ and $1<p<\infty$.  
It follows from \cite[Theorem 1.1]{laeng} that 
\begin{equation}\label{3.15}
\int_E\big|H(\chi_E)(x)\big|^pdx=\Big(2-\frac{1}{2^{p-1}}\Big)\frac{\mu(E)}{\pi^p}\zeta(p)\Gamma(p+1),
\end{equation}
where $\zeta(p)$ is the value of the Riemann zeta function $\zeta(s)$ at $s=p$.
Given any Borel  set $A\subseteq(-1,1)$ note, 
according to the definition of $H$ in \cite[(1.1)]{laeng}, that 
$$
-\big(H(\chi_A)\big)(t)=T(\chi_A)(t),\quad t\in A.
$$ 
Hence, from \eqref{3.15} we have, for each $1<p<\infty$, that 
\begin{align}\label{3.16}\nonumber 
\big\|T(\chi_A)\|_p & =
\Big(\int_{-1}^1\big|T(\chi_A)(x)\big|^pdx\Big)^{1/p}
\\ & \ge \nonumber
\Big(\int_A\big|T(\chi_A)(x)\big|^pdx\Big)^{1/p}
\\ & =
\Big(\int_A\big|H(\chi_A)(x)\big|^pdx\Big)^{1/p}
\\ & =\nonumber
\left(\Big(2-\frac{1}{2^{p-1}}\Big)\frac{\mu(A)}{\pi^p}
\zeta(p)\Gamma(p+1)\right)^{1/p}.
\end{align}
Suppose furthermore that  $\mu(A)>0$. Set $f=T(\chi_A)\in \lexp$.
Combining  \eqref{3.14} and \eqref{3.16} yields 
\begin{align*}
\big\|T(\chi_A)\big\|_{\lexp}&
\ge 
c \sup_{1\le p<\infty}\frac{\|T(\chi_A)\|_{L^p}}{p}
\\ &
\ge c \sup_{1\le p<\infty}\frac{1}{p}\left(\Big(2-\frac{1}{2^{p-1}}\Big)\frac{\mu(A)}{\pi^p}
\zeta(p)\Gamma(p+1)\right)^{1/p}	
\\ &
\ge \frac{c}{\pi}  \sup_{1\le p<\infty}\frac{1}{p}\Big(\mu(A)\Gamma(p+1)\Big)^{1/p}	
\\ &
\ge \frac{c}{\pi} \lim_{n\to\infty} \mu(A)^{1/n}\frac{(\Gamma(n+1))^{1/n}}{n}
\\ & = \frac{c}{\pi} \lim_{n\to\infty} \sqrt[n]{\frac{n!}{n^n}}
= \frac{c}{\pi }\frac{1}{ e}
> \frac{1}{e^2\pi}.
\end{align*}
\end{proof}


\begin{remark}\label{r-3.12}
Theorem \ref{t-3.11}  implies that $\mexp$ is \textit{not} $\sigma$-additive.
Actually, $\mexp$ is not even \textit{strongly additive}. Recall that 
a Banach-space-valued, finitely 
additive set function $\nu$ defined on an algebra $\mathcal F$ of  subsets of a 
non-empty set is strongly additive if,   for every sequence  
$(E_n)_{n=1}^\infty$ of pairwise disjoint sets in  $\mathcal F$, the series 
$\sum_{n=1}^\infty \nu(E_n)$ is  norm-convergent; see  \cite[Definition I.1.14]{diestel}, 
and, for equivalent conditions,  \cite[Proposition I.1.17]{diestel}.
Clearly every $\sigma$-additive vector measure is necessarily strongly additive.
\end{remark}


\section{Integral representation of  $T$}


In this section we study  operator
theoretic characteristics of $\tlog\colon\logl\to L^1$
and establish the integral representation given in  \eqref{1.1}.
The tool for such a representation is the space 
$L^1(\mlog)$ of (equivalence classes of) $\mcB$-measurable 
functions $f\colon (-1,1)\to\C$ which are integrable with respect to 
the vector measure $\mlog$.


A $\mcB$-measurable function $f\colon (-1,1)\to\C$ is said to be $\mlog$-\textit{integrable} 
if it is $\langle \mlog,\, g \rangle$-integrable for every $g \in L^\infty$ and if, 
given any $A \in \mcB$, there is  a function $\int_A f d\mlog $ belonging to  $L^1$ 
(necessarily unique) satisfying
$$
\Big\langle \int_A f d\mlog, \,g\Big\rangle = \int_A f\,d\langle \mlog,\, g 
\rangle,\quad \text{ for every    }  g \in (L^1)^*=L^\infty.
$$
In this case, $\int_Af\,d\mlog\in L^1$ is called the integral of $f$ over $A$ with 
respect to $\mlog$. 
The resulting $L^1$-valued set function
$$
A\longmapsto\int_Af\,d\mlog\in L^1,\quad A\in\mcB,
$$
is called the indefinite integral of $f$ with respect to $\mlog$; it is again 
$\sigma$-additive 
via the Orlicz-Pettis theorem, \cite[Corollary 1.4.4]{diestel}.
Every $\mcB$-simple function $\varphi=\sum_{j=1}^na_j\chi_{A_j}$ with
$a_j\in\C$ and $A_j\in\mcB$ for $j=1,\dots,n$ and $n\in\N$ is
$\mlog$-integrable and $\int_A\varphi\,d\mlog
=\sum_{j=1}^na_j \mlog(A_j\cap A)$ for $A\in\mcB$.
For the above definitions and facts
see  \cite[Section 3.1]{okada-ricker-sanchez}, for example.

By $L^1(\mlog)$ we denote the vector space of all $\mlog$-integrable 
functions.  Proposition \ref{p-3.1}(ii) implies that the indefinite 
integral of a function $f\in L^1(\mlog)$ is zero if an only if $f=0$ 
off an $\mlog$-null set if and only if
$f=0$ pointwise $\mu$-a.e. Accordingly those functions in $L^1(\mlog)$
which equal 0 $\mu$-a.e.\ are identified with the zero function. So,
$L^1(\mlog)$ is a linear subspace of $L^0$. Moreover, the space $L^1(\mlog)$
is a B.f.s.\ over $\inter$ with respect to the lattice norm
$$
\|f\|_{\lmlog}:=\sup\Big\{\int_{-1}^1|f|\,d|\langle\mlog,g\rangle|: 
g\in L^\infty \text{ with } \|g\|_\infty\le1\Big\}, \quad f\in\lmlog.
$$
For this and the fact that $\lmlog$ has $\sigma$-order
continuous norm, we refer to \cite[Theorem 1]{curbera-MA}. The denseness 
of $\text{sim }\mcB$ in $\lmlog$ is a consequence of the Lebesgue Dominated Convergence
theorem for vector measures (due to D.~R.~Lewis); see, for example, 
\cite[Theorem 3.7(ii)]{okada-ricker-sanchez}.

The integration operator $I_{\mlog}\colon \lmlog\to L^1$ is the linear map defined by
\begin{equation}\label{4.1}
I_{\mlog}(f):=\int_{(-1,1)}f\,d\mlog,\quad f\in\lmlog.
\end{equation}
It is a bounded operator as can be seen from the inequality
$$
\|I_{\mlog}(f)\|_{L^1}\le \|f\|_{\lmlog},\quad f\in\lmlog.
$$


Given any r.i.\  space $X$ over $(-1,1)$ with non-trivial Boyd indices, 
the associated  vector measure $m_X$
satisfies $L^1(m_X) = X_a$ (as  
vector sublattices of $L^0$ and with equivalent lattice norms), where
$X_a$ denotes  the absolutely continuous part of $X$, \cite[Definition I.3.1]{bennett-sharpley}.  
Moreover, $T_X$ equals the integration operator $ I_{m_X}$.  For this, 
see \cite[Theorem 3.10 and Corollary 3.11(vi)]{curbera-okada-ricker-qm} and 
\cite[Theorem]{curbera-okada-ricker-mh}.
The crucial point for this note is the fact that both $\logl$ and $L^1$ have 
absolutely continuous norm and that $\tlog\colon\logl\to L^1$ is
a bounded operator (cf.\ Lemma \ref{l-2.1}) which satisfies 
$\mlog(A)=\tlog(\chi_A)$, for $A\in\mcB$.
That is, $\mlog$ is the vector measure associated to the operator $\tlog$.


\begin{theorem}\label{t-4.1} 
The identity $L^1(\mlog) = \logl$ holds as vector sublattices of $L^0$ and
with equivalent lattice norms.   In particular,
$L^1(\mlog)$ and $\logl$ are isomorphic B.f.s.' over $(-1,1)$.  
Moreover, the integration operator $I_{\mlog}= \tlog$. 
\end{theorem}

\begin{proof}
To verify that $\logl\subseteq L^1(\mlog)$, let $f\in\logl$. 
Given $g\in L^\infty=(L^1)^*$, we have from Lemma \ref{l-2.2} that
$fT(g)\in L^1$, which implies that $f\in L^1(\langle \mlog,g\rangle)$
via \eqref{3.3}. Now, given $A\in\mcB$, since $f\chi_A\in \logl$,
Parseval's formula \eqref{2.3} with $f\chi_A$ in place of $f$ gives
$$
\int_Af\,d\langle \mlog,g\rangle= -  \int_{-1}^1f\chi_AT(g)\,d\mu
=  \int_{-1}^1gT(f\chi_A)\,d\mu=
\Big\langle \tlog(f\chi_A),g\Big\rangle.
$$
This implies that $f\in L^1(\mlog) $ and $\int_Af\,d\mlog=\tlog(f\chi_A)$
for $A\in\mcB$. In particular,  
$I_{\mlog}(f)=\tlog(f)$; see \eqref{4.1}. So, $\logl\subseteq L^1(\mlog)$.
Moreover, it has also been shown that $I_{\mlog}$ is a bounded linear 
extension of $\tlog$ to $L^1(\mlog) $.

On the other hand, $\tlog$ does not admit a bounded linear extension to
any strictly larger B.f.s.\ over $(-1,1)$, \cite[Theorem 5.6]{curbera-okada-ricker-log},
which implies that both $L^1(\mlog)=\logl$ and $I_{\mlog}=\tlog$.
That the B.f.s.' $L^1(\mlog)$ and $\logl$ have equivalent lattice norms
is a consequence of the closed graph theorem, for example.
\end{proof}


\begin{corollary}\label{c-4.2}
Let $f\in L^0$. Then $f$ belongs to $\logl$ if and only if $fT(g)\in L^1$ for 
every $g\in L^\infty$.
\end{corollary}

\begin{proof}
$\Rightarrow$) This is part of Lemma \ref{l-2.2}.

$\Leftarrow$) According to \eqref{3.3} the condition  $fT(g)\in L^1$ means that
$f\in L^1(\langle\mlog,g\rangle)$ for every $g\in L^\infty=(L^1)^*$.
In other words, $f$ is scalarly $\mlog$-integrable. Since the codomain space $L^1$
of $\mlog$ 
does not contain an isomorphic copy of $c_0$
(as it is weakly sequentially complete), every scalarly $\mlog$-integrable function is necessarily $\mlog$-integrable, \cite[Theorem 5.1]{lewis}. 
Hence, $f\in\lmlog=\logl$ by Theorem \ref{t-4.1}.
\end{proof}


\begin{remark}\label{r-4-3}
(i) The optimal domain $[T,L^1] $, for $T$ taking values in $L^1$, has 
been specified in \cite[(5.1) and (5.2)]{curbera-okada-ricker-log}.  
Due to \cite[Lemma 5.4]{curbera-okada-ricker-log} we can conclude  
that $[T,L^1]= L^1(\mlog)$ as identical B.f.s.' over $(-1,1)$.

(ii) From  \cite[Proposition 5.1(i) and Theorem 5.6]{curbera-okada-ricker-log}
it follows that a function $f\in L^1$ belongs to $\logl$ if and only if
$fT(g)\in L^1$ for every $g\in L^\infty$. Corollary \ref{c-4.2} above 
also presents this fact (even more general,  for $f\in L^0$), possible by using  the vector measure $\mlog$.

(iii) It follows from Theorem \ref{t-4.1} that $L^1(\mlog)$ is a \textit{r.i.\ space}. 
This is not true in general for the $L^1$-space of the vector measure associated with a kernel operator, \cite[Example 5.15(b)]{curbera-ricker}.
\end{remark}


The identification of  $\tlog\colon\logl\to L^1$ with 
the integration operator $I_{\mlog}\colon\lmlog\to L^1$, established in Theorem \ref{t-4.1}, 
allows us to prove the following operator theoretic results concerning $\tlog$.


\begin{theorem}\label{t-4.4} 
The   finite Hilbert transform  $\tlog\colon\logl\to L^1$ 
is not  order bounded.
\end{theorem}

\begin{proof}
Since $\{\chi_B:B\in\mcB\}$ is an order bounded subset of $\logl$ but
$\mlog(\mcB)=\{\tlog(\chi_B):B\in\mcB\}$ is \textit{not} an order bounded 
subset of $L^1$
(cf.\ Corollary \ref{c-3.6} with $A=(-1,1)$), it follows that $\tlog\colon\logl\to L^1$ 
is not an order bounded operator.
\end{proof}


We now provide two alternative arguments for the non-compactness 
of $\tlog$ to that given in 
\cite[Proposition 4.15]{curbera-okada-ricker-log}.

\begin{corollary}\label{c-4.5}
The   finite Hilbert transform $\tlog\colon\logl\to L^1$ is not compact.
\end{corollary}

\begin{proof}
It suffices to establish that the integration operator 
$I_{\mlog}\colon\lmlog\to L^1$ is not compact.
But, according to Theorems 1 and 4 of  \cite{ORR}, a Banach-space-valued
measure must have finite variation whenever its 
associated integration operator is compact.
This fact and Theorem \ref{t-3.4} prove the result.

An alternative proof follows from Theorem 
\ref{t-3.10} and the fact that $\mlog(\mcB)=\tlog(\{\chi_A:A\in\mcB\})$ 
with $\{\chi_A:A\in\mcB\}$ 
a bounded subset of $\logl$.
\end{proof}


A linear operator between Banach spaces is called 
\textit{completely continuous} if it maps every weakly compact 
set into a norm-compact set.
It is a concept that lies between that of a bounded operator 
and that of a compact operator.


\begin{theorem}\label{t-4.6}
The   finite Hilbert transform $\tlog\colon\logl\to L^1$ is not completely continuous.
\end{theorem}

\begin{proof}
It suffices to establish that the integration operator 
$I_{\mlog}\colon\lmlog\to L^1$ is not completely continuous.
But, this  holds via Theorem \ref{t-3.10} by applying  the general fact that
if the range of a Banach-space-valued vector measure  is not relatively compact, 
then its associated integration operator is not completely continuous; 
see, for example, \cite[p.153]{okada-ricker-sanchez}.  
\end{proof}

Due to Corollary \ref{c-4.5}, the integration operator $I_{\mlog}$ is not compact. 
This raises the question as to whether or not $I_{\mlog}$ is at least weakly compact.   
The following result provides the answer.


\begin{theorem}\label{t-4.7} 
The   finite Hilbert transform $\tlog\colon\logl\to L^1$
is not weakly compact.
Consequently, neither is  the integration operator $I_{\mlog}\colon   L^1(\mlog) \to L^1$.
\end{theorem}

\begin{proof}
We shall show  that the adjoint operator 
$(T_{\mathrm{log}})^*\colon L^\infty\to \lexp$ of $T_{\mathrm{log}}$  
is not weakly compact.  By \eqref{2.5}, we have 
$(T_{\mathrm{log}})^* = -T_\infty$.
According to Remark \ref{r-3.12}  
the bounded finitely additive set function
$\mexp\colon  A\mapsto T_\infty(\chi_A) \in \lexp$ for $  A \in  \mathcal B$ 
is not strongly additive. Hence,   
its range $\mexp(\mathcal B)$ is not relatively weakly 
compact in $\lexp$; see  \cite[Corollary I.5.3]{diestel}.   
Consequently, the image of the unit ball $B_{L^\infty}$ 
under $T_\infty$ is not relatively weakly compact in $\lexp$ 
because $\mexp(\mathcal B) \subseteq T_\infty(B_{L^\infty})$.
Thus,  $(T_{\mathrm{log}})^* = -T_\infty$ is not a weakly compact operator.  
We can conclude from Gantmacher's theorem, \cite[Theorem VI.4.8]{DS}, 
that also $T_{\mathrm{log}}$ is not weakly compact.

Consequently, 
the integration operator $I_{\mlog}\colon   L^1(\mlog) \to L^1$
is not weakly compact either.
\end{proof}



\end{document}